
\input amstex
\documentstyle{amsppt}

\hcorrection{19mm}

\nologo
\NoBlackBoxes


\topmatter 

\title Thin position and essential planar surfaces
\endtitle
\author  Ying-Qing Wu$^1$
\endauthor
\leftheadtext{Ying-Qing Wu}
\rightheadtext{Thin position and planar essential surfaces}
\address Department of Mathematics, University of Iowa, Iowa City, IA
52242
\endaddress
\email  wu\@math.uiowa.edu
\endemail
\keywords Thin position, knots and links, essential planar surfaces
\endkeywords
\subjclass  Primary 57M25
\endsubjclass

\thanks  $^1$ Partially supported by NSF grant \#DMS 0203394
\endthanks

\abstract Abby Thompson proved that if a link $K$ is in thin position
but not in bridge position then the knot complement contains an
essential meridional planar surface, and she asked whether some thin
level surface must be essential.  This note is to give a positive
answer to this question, showing that the if a link is in thin
position but not bridge position then a thinnest level surface is
essential.  A theorem of Rieck and Sedgwick follows as a consequence,
which says that thin position of a connected sum of small knots comes
in the obvious way.  \endabstract

\endtopmatter

\document

\define\proof{\demo{Proof}}
\define\endproof{\qed \enddemo}
\define\a{\alpha}
\redefine\b{\beta}

\redefine\e{\epsilon}
\redefine\bdd{\partial}
\define\Int{\text{\rm Int}}
\input epsf.tex

The concept of thin position was introduced by David Gabai in [G], and
has been used successfully in attacking some very difficult problems,
see for example [G, GL, ST, T1].  In [T2] Abby Thompson proved that if
a knot $K$ is in thin position but not in bridge position then some
thin level surface can be compressed to produce an essential planar
surface in the complement of $K$ with meridional boundary slope; in
particular, by a theorem of Culler, Gordon, Luecke and Shalen [CGLS]
this implies that the knot is large in the sense that its complement
contains some closed essential surfaces.  This was further explored by
Heath and Kobayashi [HK], showing that certain thin level surface of a
nonsplit link $L$ in thin position but not bridge position can be
compressed to give some natural tangle decomposition of $L$, and the
decomposing spheres then give rise to essential meridional planar
surfaces in the link complement.  In both [T2] and [HK] the essential
planar surfaces come from compression of a thin level surface.  There
are examples in [HK] showing that some essential meridional planar
surfaces are not level surface of a knot in thin position.  This leads
to a question raised by Thompson [T3], which asks whether 
a link $L$ in thin position but not in bridge position has a level
surface which is essential.  The purpose of this paper is to give a
positive solution to this problem.

\proclaim{Theorem 1} If a link $L$ in $S^3$ is in thin position but
not in bridge position, then a thinnest level surface $Q$ of $L$ is an
essential surface in $S^3 - \Int N(L)$.  \endproclaim

We give some definitions.  Consider $S^3$ as $\Bbb R^3 \cup
\{\infty\}$, and let $\rho$ be the height function $\rho: \Bbb R^3 \to
\Bbb R$ defined by $\rho(x,y,z) = z$.  For each $t\in \Bbb R$, let
$P(t) = \rho^{-1}(t)$ be the horizontal plane in $\Bbb R^3$ at height
$t$.  When $t$ is not a critical level of $\rho$, define $Q(t)$ to be
the punctured sphere $(P(t) \cup \{\infty\}) - \Int N(L)$, called the
{\it level surface\/} at level $t$, where $N(L)$ is a small regular
neighborhood of $L$ intersecting $P$ in meridional disks.  If $I$ is
an interval on $\Bbb R$, denote by $Z_I = \Bbb R^2 \times I$ the set
of points in $\Bbb R^3$ whose $z$-coordinate is in $I$.

Let $L$ be a link in $\Bbb R^3$ such that the restriction of $\rho$ to
$L$ is a Morse function, and let $a_0, ..., a_n$ be the critical points
of $L$, labeled so that the corresponding critical values $t_i =
\rho(a_i)$ satisfy $t_{i-1} < t_i$ for all $i$.  Let $s_i \in
(t_{i-1}, t_i)$.  Thus $P_i = P(s_i)$ is a plane between $a_{i-1}$ and
$a_i$, called a level plane corresponding to the critical point
$a_{i-1}$.  The {\it width of $P(t)$\/} (with respect to $L$) is
defined as $w(P(t)) = |P(t) \cap L|$, where $|A|$ denotes the number
of elements in $A$.  The {\it width of $L$\/} is $w(L) = \Sigma_1^n
w(P_i)$.  A link $L$ is in {\it thin position\/} if $w(L)$ is minimal
up to isotopy of $L$.

A plane $P(t)$ with $t_{i-1} < t < t_i$ is called a {\it thin level
plane\/} of $L$, and $t$ a thin level, if $a_{i-1}$ is a local maximum
and $a_i$ a local minimum.  Similarly, $t\in (t_{i-1}, t_i)$ is a {\it
thick level\/} and $P(t)$ a thick level plane if $a_{i-1}$ is a local
minimum and $a_i$ a local maximum.  A thin level plane $P(t)$ is a
{\it thinnest level plane\/} and the corresponding planar surface
$Q(t) = (P(t)\cup \{\infty\}) - \Int N(L)$ a {\it thinnest level
surface\/} if $w(P(t))$ is minimal among all thin planes.  The link
$L$ is {\it in bridge position\/} if it has no thin level.

\medskip
\demo{Proof of Theorem 1} Let $P$ be a thinnest level plane, and $Q$
the corresponding planar surface defined above.  Moving $L$ up or down
if necessary we may assume without loss of generality 
that $P = P(0)$, i.e., it is the $xy$-plane in $\Bbb R^3$.  Our goal
is to show that $Q$ is an essential surface in $E(L)$.  Recall that a
properly embedded compact orientable surface in a 3-manifold $M$ is
{\it essential\/} if (i) it is incompressible (in particular it is not
a 2-sphere bounding a 3-ball), and (ii) it is not boundary parallel.
Since $P$ is a thin level plane, it is easy to show that $Q$ is never
a 2-sphere bounding a 3-ball, and it is not boundary parallel.
Therefore we need only show that it is not a compressible punctured
sphere.

Assume to the contrary that $Q$ is compressible and let $D$ be a
compressing disk.  Without loss of generality we may assume that $D$
is in the upper half space $\Bbb R^3_+ = Z_{[0,\infty)}$.  We will show
below that either $L$ is not in thin position or $P$ is not a thinnest
level plane, which will contradict the assumption and complete the
proof of the theorem.

The $xz$-coordinate plane cuts $\Bbb R^3_+$ into $Y_-$ and $Y_+$,
where $Y_-$ is the left half space of $\Bbb R^3_+$, consisting of
points in $\Bbb R^3_+$ with negative $y$-coordinate, and $Y_+$ the
right half space of $\Bbb R^3_+$.  Let $D'$ be the disk on $P$ with
$\bdd D = \bdd D'$, and let $B$ be the 3-ball bounded by $D \cup D'$.
Define $\a = L \cap B$, and $\b = L \cap (\Bbb R^3_+ - B)$.  Since $D$
is a compressing disk of $Q$, both $\a$ and $\b$ are nonempty.

An isotopy $\phi_t$ of $\Bbb R^3$ is called an {\it h-isotopy\/} if
$\phi_0 = id$, and $\phi_1$ is a level preserving map on $L$, i.e.,
$\rho \circ \phi_i = \rho$ on $L$.  Note that $\phi_t$ does not have
to be level preserving on $L$ when $t \neq 0, 1$.  Since $\phi_1$ is
level preserving, the link $\phi_1(L)$ has the same width as $L$.

\proclaim{Lemma 2}  
Up to h-isotopy we may assume that $\a \subset Y_-$ and $\b \subset Y_+$.
\endproclaim

\proof By a level preserving isotopy which shifts the whole link $L$
to the left we may assume that the 3-ball $B$ lies in $Y_-$.  Let
$I=[0,\e]$ be an interval containing no critical value of $\rho$, so
$L \cap Z_I$ can be assumed to be a set of vertical arcs from $\Bbb R^2
\times 0$ to $\Bbb R^2 \times \e$.  We may also assume that $B \cap
Z_I = D' \times I$.  Let $f_t$ be an isotopy supported in $B$ which
shrinks $\a$ into $D'\times I$, i.e., $f_0 = id$ and $f_1(\a) \subset
D' \times I$.  (Note that $f_t$ is not level preserving.)  There is
now a level preserving isotopy $g_t$ of $f_1(L)$, supported outside of
$D' \times I$, moving $\b$ into $Y_+$.  Let $h_t$ be the reverse
isotopy of $f_t$, i.e., $h_t = f_{1-t}$.  Then the union of these
three isotopies is the required h-isotopy.  \endproof

The separation of $\a$ and $\b$ by the $xz$-coordinate plane allows us
to modify $\a$ by ``vertical isotopy'' without intersecting the rest
of $L$.  Let $\phi_t$ be an isotopy of the positive half of the
$z$-axis.  Then $id \times \phi_t$ is an isotopy of $Y_- = \Bbb R^2_-
\times \Bbb R_+$, which can be extended to an isotopy $f_t$ of $\Bbb
R^3$ supported in a small neighborhood of $Y_-$ in $\Bbb R^3_+$, and
hence is the identity on $L - \a$.  This $f_t$ is called a vertical
isotopy on $Y_-$ determined by $\phi_t$.  

Denote by $m_{\a}$ the maximum value of $\rho(\a)$. Let $n_{\a}$ be
the first local minimum level of $\a$, counted from top down, and
$n_{\a} = 0$ if $\a$ has no local minimum.  (It can be shown that $\a$
must have some local minima, but this is not necessary.)  Similarly
for $m_{\b}$ and $n_{\b}$.

\proclaim{Lemma 3} Either $n_{\a} > m_{\b}$ or $n_{\b} > m_{\a}$.
\endproclaim

\proof Since $P$ is a thin level plane, at least one of $\a$ or $\b$
contains some local minima, so $n_{\a}$ and $n_{\b}$ cannot both be
$0$.  Without loss of generality we may assume that $n_{\a} < n_{\b}$.
Then above the level $n_{\b}$ all the critical points of $L$ are local
maxima.  Note that lowering a local maximum of $\a$ through the level
of a local maximum of $\b$ will not change $w(L)$.  If $m_{\a} >
n_{\b}$ then deforming a local maximum of $\a$ downward to a level
just below $n_{\b}$ would reduce $w(L)$, contradicting the minimality
of $w(L)$.  \endproof

By Lemma 3, we may assume without loss of generality that $m_{\a} <
n_{\b}$.  In particular $\b$ must have some minima.  Let $r \in
(0,n_{\b})$ be such that $|P(r)\cap \b|$ is minimal among all $|P(z)
\cap \b|$, $z \in [0,n_{\b}]$.  Let $I=[a,b]$ be a maximal interval in
$[0,\infty)$ containing $r$ and having no critical value of $\b$ in
its interior.  By the definition of $r$, one can see that $b$ is a
local minimum level of $\b$ and $a$ is either $0$ or a local maximum
level of $\b$.

Let $f_t$ be a vertical isotopy supported in $Z_{[a,\infty)}$, which
pushes $\a$ downward to $\a' = f_1(\a)$ lying below the level $b$.
Let $L'$ be the presentation of $L$ obtained this way.  We will show
below that $w(L') < w(L)$, which then contradicts the assumption,
completing the proof of Theorem 1.

Since the isotopy $f_t$ is supported in $Z_{[a, \infty)}$, the
critical points of $L$ below level $a$ and the widths of the
corresponding level planes remain unchanged, so we need only calculate
the sum of the width for those level planes corresponding to critical
points of $L$ and $L'$ above level $a$.  Denote by $P_1, ..., P_k$
(resp.\ $P'_1, ..., P'_k$) the level planes of $L$ corresponding to
(i.e., lying just above) the critical points of $\a$ (resp.\ $\a'$)
above the level $a$, and $R_1, ..., R_h$ those corresponding to
critical points of $\b$ above level $a$, labeled according to their
height.  Since $\b$ has no critical values between $a$ and $b$, all
$R_j$ are above level $b$, so we have
$$|R_j \cap L| = |R_j \cap (\a \cup \b)| \geq |R_j \cap \b| = |R_j
\cap L'|.$$ Also, since the top level of $\a$ is below $n_{\b}$,
by the choice of $r$ and the fact that $|P_i \cap \a| = |P'_i \cap
\a|$ we have
$$ \align |P_i \cap (\a \cup \b)| &= |P_i \cap \a| + |P_i \cap \b| 
\geq |P_i \cap \a| + |P(r) \cap \b| \\
&= |P'_i \cap \a| + |P'_i \cap \b| = |P'_i \cap L'|. \endalign$$
It follows that $w(L) \geq w(L')$, and equality holds if and only if
it holds in all the above inequalities.  On the other hand, since $P$
is a thinnest level of $L$, the level plane just below $b$
cannot be disjoint from $\a$ because otherwise it would be a thin
level plane with width $|P(r) \cap \b| \leq |P \cap \b| < |P \cap L|$,
contradicting the choice of $P$.  Therefore $\a$ must intersect
the level plane $R_1$ lying just above level $b$,  so $|R_1
\cap L| > |R_1 \cap L'|$, and $w(L) > w(L')$, which contradicts the
assumption that $L$ is in thin position.
\endproof

Let $K_1$ and $K_2$ be knots in $S^3$.  Putting in thin position, with
$K_1$ above $K_2$, then taking the obvious connected sum, we have a
projection of $K= K_1 \# K_2$ with width $w(K_1) + w(K_2) -2$.  Hence
we have $w(K) \leq w(K_1) + w(K_2) - 2$.  It was conjectured that
$w(K) = w(K_1) + w(K_2) - 2$ for all $K = K_1 \# K_2$.  See [RS, SS]
for some work concerning this conjecture.

Recall that a knot $K$ is a {\it small knot\/} if its complement
contains no closed essential surface.  An essential surface $F$ in the
knot exterior $E(K) = S^3 - \Int N(K)$ is a {\it meridional essential
surface\/} if $\bdd F$ is a nonempty set of meridional curves on $\bdd
E(K)$.  By [CGLS], if $K$ is small then its exterior contains no
meridional essential planar surface.  The following result is due to
Rieck and Sedgwick [RS].  It proves the above conjecture for connected
sum of small knots $K_1$ and $K_2$.

\proclaim{Corollary 4 ([RS])}  Let $K_1$ and $K_2$ be nontrivial knots in
$S^3$ such that $E(K_i) = S^3 - \Int N(K_i)$ contains no meridional
essential planar surfaces.  Let $K$ be a thin position embedding of
$K_1 \# K_2$.  Then there is a level sphere $S$ in $S^3$ intersecting
$K$ in two points, decomposing $K$ into $K_1$ and $K_2$.  In
particular, $w(K_1 \# K_2) = w(K_1) + w(K_2) - 2$.  \endproclaim

\proof Let $F$ be sphere in $S^3$ which realizes the connected sum
$K_1 \# K_2$.  Then $P = F \cap E(K)$ is a meridional essential planar
surface.  Since the exterior of $K_i$ contains no meridional essential
planar surface, one can show that $P$ is the only meridional essential
planar surface in $E(K)$.

It is easy to see that $K$ cannot be in bridge position: By Schubert's
theorem $b(K) = b(K_1) + b(K_2) -1$, so a bridge position projection
of $K$ can be obtained by putting $K_i$ in bridge position, with $K_1$
above $K_2$, taking the connected sum, then raising the maxima of
$K_2$ over the minima of $K_1$.  But the last operation would increase
width, hence if $K = K_1 \# K_2$ is in bridge position then it cannot
be in thin position.  (See [RS] for an alternative proof.)

It now follows from Theorem 1 that if $K$ is in thin position then
some thinnest level surface $Q$ of $K$ is an essential planar surface
in $S^3 - \Int N(K)$.  Since meridional essential planar surface of
$K$ is unique, $Q$ is the same as the surface $P$ above up to isotopy,
hence the result follows.  
\endproof

\medskip

{\it Acknowledgment.\/} I would like to thank Abby Thompson for
raising the question, and for some interesting conversation about the
proof of Theorem 1.  Thanks also to Jennifer Schultens for some
interesting discussion about thin position and width of connected sum
of knots.

\enddocument